\documentclass{amsart}
\usepackage{amsmath}
\usepackage{amssymb}
\usepackage{latexsym}
\usepackage{amsbsy}
\usepackage[all]{xy}

\usepackage{amscd}
\usepackage[dvips]{graphicx}
\newtheorem{theorem}{Theorem}[section]
\newtheorem{lemma}[theorem]{Lemma}

\theoremstyle{definition}

\newtheorem{Prop}[theorem]{Proposition}
\newtheorem{Cor}[theorem]{Corollary}

\theoremstyle{remark}

\numberwithin{equation}{section}

\newcommand{\bbe}{\mathbb{E}}
\newcommand{\str}[1]{\langle #1\rangle}
\newcommand{\az}{\alpha}

\newcommand{\ud}{{\underline{d}}}

\newcommand\lz{\lambda}

\newcommand\bbn{{\mathbb N}}

\newcommand\bbc{{\mathbb C}}

\newcommand\ext{\mbox{Ext}\,}

\newcommand\nd{{\noindent}}

\newcommand\mc{{\mathcal{C}}}

\newcommand\mo{{\mathcal{O}}}

\newcommand\ue{{\underline{e}}}

\newcommand{\id}{1\kern -.35em 1}

\begin{document}

\title{A multiplication formula for module subcategories of Ext-symmetry}

\author{Jie Xiao}
\address{Department of Mathematics, Tsinghua University, Beijing 100084, P.R.China}
\email{jxiao@math.tsinghua.edu.cn}
\author{Fan Xu}
\address{Department of Mathematics, Tsinghua University, Beijing 100084, P.R.China}

\email{fanxu@mail.tsinghua.edu.cn}
\thanks{The research was
supported in part by NSF of China (No. 10631010)}


\date{December 30, 2007, last revised October 30, 2008}


\begin{abstract}
We define evaluation forms associated to objects in a module
subcategory of Ext-symmetry generated by finitely many simple
modules over a path algebra with relations and prove a
multiplication formula for the product of two evaluation forms. It
is analogous to a multiplication formula for the product of two
evaluation forms associated to modules over a preprojective
algebra given by Geiss, Leclerc and Schr\"oer in \cite{GLS2006}.
\end{abstract}

\maketitle

\section*{Introduction}
Let $\Lambda$ be the preprojective algebra associated to a
connected quiver without loops (see e.g. \cite{Ringel}) and
$\mathrm{mod}(\Lambda)$ be the category of finite-dimensional
nilpotent left $\Lambda$-modules. We denote by $\Lambda_{\ue}$ the
variety of finite-dimensional nilpotent left $\Lambda$-modules
with dimension vector $\ue.$ For any $x\in \Lambda_{\ue},$ there
is an evaluation form $\delta_x$ associated to $x$ satisfying that
there is a finite subset $R(\ue)$ of $\Lambda_{\ue}$ such that
$\Lambda_{\ue}=\bigsqcup_{x\in R(\ue)}\str{x}$ where
$\str{x}:=\{y\in\Lambda_{\ue}\mid \delta_{x}=\delta_{y} \}$
\cite[Section 1.2]{GLS2006}. Inspired by the Caldero-Keller
cluster multiplication theorem for finite type \cite{CK2005},
Geiss, Leclerc and Schr\"oer \cite{GLS2006} proved a
multiplication formula (the Geiss-Leclerc-Schr\"oer multiplication
formula) as follows:
$$
\chi(\mathbb{P}\ext^1_\Lambda(x',x''))\, \delta_{x'\oplus x''} =
\sum_{x\in R(\ue)}
\left(\chi(\mathbb{P}\ext^1_\Lambda(x',x'')_{\str{x}})+
\chi(\mathbb{P}\ext^1_\Lambda(x'',x')_{\str{x}})\right)\delta_x,
$$
where  $x'\in \Lambda_{\ue'},$ $x''\in \Lambda_{\ue''}$,
$\ue=\ue'+\ue'',$ $\mathbb{P}\ext^1_\Lambda(x',x'')_{\str{x}}$ is
the constructible subset of $\mathbb{P}\ext^1_\Lambda(x',x'')$
with the middle terms belonging to $\str{x},$
$\mathbb{P}\ext^1_\Lambda(x'',x')_{\str{x}}$ is defined similarly.

The proof of the formula heavily depends on the fact that the
category $\mathrm{mod}(\Lambda)$ is of Ext-symmetry. A category
$\mathcal{C}$ is of Ext-symmetry if there is a bifunctorial
isomorphism: $\mathrm{Ext}^1_{\mc}(M,N)\cong
\mathrm{DExt}^1_{\mc}(N,M)$ for any objects $M, N\in \mc$.

Let $Q$ be a finite quiver and $A$ be a quotient algebra $\bbc
Q/\mathcal{I}$ by an ideal $\mathcal{I}$. We denote by
$\mathrm{mod}(A)$ be the category of finite dimensional left
$A$-modules. We call $A$ an algebra of Ext-symmetry if
$\mathrm{mod}(A)$ is of Ext-symmetry.  It is proved that
preprojective algebras and deformed preprojective algebras are of
Ext-symmetry (see \cite[Theorem 3]{GLS2006} and Section 3 in this
paper).

In this paper, we focus on the module subcategories of
Ext-symmetry of $\mathrm{mod}(A)$. Let $\mathcal{S}=\{S_1,\cdots,
S_n\}$ be a finite subset of finite-dimensional simple
$A$-modules. We denote by $\mc(\mathcal{S})$ the full subcategory
of $\mathrm{mod}(A)$ consisting of modules $M$ satisfying that the
isomorphism classes of the composition factors of $M$ belong to
$\mathcal{S}$. We associate to modules in $\mc(\mathcal{S})$ some
evaluation forms and prove that if $\mc(\mathcal{S})$ is of
Ext-symmetry, then the product of two evaluation forms satisfies
an identity (Theorem 2.3). The identity is analogous to the
Geiss-Leclerc-Schr\"oer multiplication formula. There are no known
examples of algebras of Ext-symmetry, apart from preprojective and
deformed preprojective algebras (see Section 3), and it is an open
question whether further examples exist. However, other examples
of module subcategories of Ext-symmetry can be easily constructed,
and we give an example in Section 3.

\section{The product of two evaluation forms}\label{basic}
\subsection{Module varieties} Let $Q=(Q_0,Q_1,s,t)$ be a finite connected quiver where $Q_0$ and $Q_1$ are the sets of vertices and
arrows, respectively, and $s,t: Q_1\rightarrow Q_0$ are maps such
that any arrow $\az$ starts at $s(\az)$ and terminates at
$t(\az).$ The space spanned by all paths of nonzero length is a
graded ideal of $\bbc Q$ and we will denote it by $\mathcal{J}.$ A
relation for $Q$ is a linear combination
$\sum_{i=1}^r\lambda_ip_i$ where $\lambda_i\in \bbc$ and the $p_i$
are paths with $s(p_i)=s(p_j)$ and $t(p_i)=t(p_j)$ for any $1\leq
i, j\leq r.$  Here if $p_i$ is a vertex in $Q_0$, then
$s(p_i)=t(p_i)=p_i$. Let $A=\bbc Q/\mathcal{I}$ where
$\mathcal{I}$ is an ideal generated by a finite set of relations.
We don't assume that $\mathcal{I}$ is admissible, i.e.
$\mathcal{I}\subset \mathcal{J}^2.$

A dimension vector for $A$ is a map $\ud: Q_0\rightarrow \bbn$. We
write $d_i$ instead of $d(i)$ for any $i\in Q_0.$ For any
dimension vector
 $\ud=(d_i)_{i\in Q_0},$ we consider the affine space over $\bbc$
$$\bbe_{\ud}(Q)=\bigoplus_{\az\in Q_1}\hom_{\bbc}(\bbc^{d_{s(\az)}},\bbc^{d_{t(\az)}}).$$
Any element $x=(x_{\az})_{\az\in Q_1}$ in $\bbe_{\ud}(Q)$ defines
a representation $(\bbc^{\ud}, x)$ where
$\bbc^{\ud}=\bigoplus_{i\in Q_0}\bbc^{d_i}$. For any
$x=(x_{\az})_{\az\in Q_1}\in\bbe_{\ud}(Q)$ and any path
$p=\az_1\az_2\cdots\az_m$ in $Q$, we set
$x_{p}=x_{\az_1}x_{\az_2}\cdots x_{\az_m}.$ Then $x$ satisfies a
relation
 $\sum_{i=1}^{r}\lz_{i}p_i$ if $\sum_{i=1}^{r}\lz_i x_{p_i}=0.$ Here if $p_i$ is a vertex in $Q_0$, then $x_{p_i}$ is the identity matrix.
Let $R$ be a finite set of relations generating the ideal
$\mathcal{I}$. Then we denote by  $\bbe_{\ud}(A)$ be the closed
subvariety of $\bbe_{\ud}(Q)$ which consists of elements
satisfying all relations in $R.$

Let $\mathcal{S}=\{S_1,\cdots, S_n\}$ be a finite subset of
finite-dimensional simple $A$-modules and $\mc(\mathcal{S})$ be a
 module subcategory of Ext-symmetry of $\mathrm{mod}(A)$. We denote by
$A_{\ud}(\mathcal{S})$ the constructible subset of $\bbe_{\ud}(A)$
consisting of modules in $\mathcal{C}(\mathcal{S})$. In the
sequel, we will fix the finite set $\mathcal{S}$ and write
$A_{\ud}$ instead of $A_{\ud}(\mathcal{S})$. The algebraic group
$G_{\ud}:=G_{\ud}(Q)=\prod_{i\in Q_0}\mathrm{GL}_{d_i}(\bbc)$ acts
on
 $\bbe_{\ud}(Q)$ by $(x_{\az})^{g}_{\az\in Q_1}=(g_{t(\az)}x_{\az}g_{s(\az)}^{-1})_{\az\in Q_1}$ for $g\in G_{\ud}$ and
  $(x_{\az})_{\az\in Q_1}\in\bbe_{\ud}(Q).$
It naturally induces the action of $G_{\ud}$ on
$A_{\ud}(\mathcal{S}).$ The orbit space is denoted by
$\overline{A}_{\ud}(\mathcal{S})$. A constructible function over
$\bbe_{\ud}(A)$  is a function $f: \bbe_{\ud}(A)\rightarrow \bbc$
such that $f(\bbe_{\ud}(A))$ is a finite subset of $\bbc$ and
$f^{-1}(c)$ is a constructible subset of $\bbe_{\ud}(A)$ for any
$c\in Q.$

Throughout this paper, we always assume that $\mc(\mathcal{S})$ is
of Ext-symmetry and constructible functions over $\bbe_{\ud}(A)$
are $G_{\ud}$-invariant for any dimension vector $\ud$ unless
particularly stated.

\subsection{Euler characteristics}
Let $\chi$ denote the Euler characteristic in compactly-supported
cohomology. Let $X$ be a complex algebraic variety and $\mo$ a
constructible subset as the disjoint union of finitely many
locally closed subsets $X_i$ for $i=1,\cdots,m.$ Define
$\chi(\mo)=\sum_{i=1}^m\chi(X_i).$ We note that it is
well-defined. The following properties will be applied to compute
Euler characteristics.
\begin{Prop}[\cite{Riedtmann} and \cite{Joyce}]\label{Euler} Let $X,Y$ be algebraic varieties over $\mathbb{C}.$
Then
\begin{enumerate}
    \item  If an algebraic variety $X$ is the disjoint union of
finitely many constructible sets $X_1,\cdots,X_r$, then
$$\chi(X)=\sum_{i=1}^{r}{\chi(X_i)}.$$
    \item  If $\varphi:X\longrightarrow Y$ is a morphism
with the property that all fibers have the same Euler
characteristic $\chi$, then $\chi(X)=\chi\cdot \chi(Y).$ In
particular, if $\varphi$ is a locally trivial fibration in the
analytic topology with fibre $F,$ then $\chi(Z)=\chi(F)\cdot
\chi(Y).$
    \item $\chi(\bbc^n)=1$ and $\chi(\mathbb{P}^n)=n+1$ for all $n\geq
    0.$
\end{enumerate}
\end{Prop}
We recall the {\it pushforward} functor from the category of
algebraic varieties over $\mathbb{C}$ and the category of
$\bbc$-vector spaces (see \cite{Macpherson} and \cite{Joyce}). Let
$\phi: X\rightarrow Y$ be a morphism of varieties. Write $M(X)$
for the $\bbc$-vector space of constructible functions on $X$. For
$f\in M(X)$ and $y\in Y,$ define
$$
\phi_{*}(f)(y)=\sum_{c\neq 0}c\chi(f^{-1}(c)\cap \phi^{-1}(y)).
$$

\begin{theorem}[\cite{Dimca},\cite{Joyce}]\label{Joyce}
Let $X,Y$ and $Z$ be algebraic varieties over $\mathbb{C},$ $\phi:
X\rightarrow Y$ and $\psi: Y\rightarrow Z$ be morphisms of
varieties, and $f\in M(X).$ Then $\phi_{*}(f)$ is constructible,
$\phi_{*}: M(X)\rightarrow M(Y)$ is a $\bbc$-linear map and
$(\psi\circ \phi)_{*}=(\psi)_{*}\circ (\phi)_{*}$ as $\bbc$-linear
maps from $M(X)$ to $M(Z).$
\end{theorem}

\subsection{The actions of $\bbc^*$ on the
extensions and flags}\label{C*}
 Let
$A=\bbc Q/\str{R}$ be an algebra as in Section 1.1. For any
$A$-modules $X,Y,$ let $D(X,Y)$ be the vector space over $\bbc$ of
all tuples $d=(d(\alpha))_{\alpha\in Q_1}$ such that linear maps
$d(\alpha)\in \mathrm{Hom}_{\bbc}(X_{s(\alpha)},Y_{t(\alpha)})$
and the matrices $
L(d)_{\alpha}=\left(%
\begin{array}{cc}
  Y_{\alpha} & d(\alpha) \\
  0 & X_{\alpha} \\
\end{array}%
\right) $ satisfy the relations in $R.$ Define $\pi:
D(X,Y)\rightarrow \mathrm{Ext}^1(X,Y)$ by sending $d$ to the
equivalence class of the following short exact sequence
$$
\xymatrix{\varepsilon:\quad 0\ar[r]& Y\ar[rr]^{\left(%
\begin{array}{c}
  1 \\
  0 \\
\end{array}%
\right)}&&L(d)\ar[rr]^{\left(%
\begin{array}{cc}
  0 & 1 \\
\end{array}%
\right)}&&X\ar[r]&0},
$$
where, as a vector space, $L(d)=(L(d)_{\alpha})_{\alpha\in Q_1}$
is the direct sum of $Y$ and $X$. The direct computation shows
that $\mathrm{Ker}\pi$ is the subspace of $D(X,Y)$ consisting of
all tuples $d=(d(\alpha))_{\alpha\in Q_1}$ such that there exist
$(\phi_i)_{i\in Q_0}\in \bigoplus_{i\in
Q_0}\mathrm{Hom}_{\bbc}(X_i, Y_i)$ satisfying
$d(\alpha)=\phi_{t(\alpha)}X_\alpha-Y_{\alpha}\phi_{s(\alpha)}$
for all $\alpha\in Q_1$ (see \cite[Section 5.1]{GLS2006} for a
similar discussion).

Fix a vector space decomposition $D(X,Y)=\mathrm{Ker}\pi\oplus
E(X,Y).$ We can identify $\mathrm{Ext}_A^{1}(X,Y)$ with $E(X,Y)$
\cite{Riedtmann}\cite{DXX}\cite{GLS2006}. Let
$\mathrm{Ext}^{1}_A(X,Y)_{L}$ be the subset of
$\mathrm{Ext}^{1}_A(X,Y)$ with the middle term isomorphic to $L,$
then $\mathrm{Ext}^{1}(X,Y)_{L}$ can be viewed as a constructible
subset of $\mathrm{Ext}_A^1(X,Y)$ by the identification between
$\mathrm{Ext}_A^1(X,Y)$ and $E(X,Y).$ There is a natural
$\bbc^{*}$-action on $E(X,Y)\setminus\{0\}$ by $t.d=(td(\alpha))$
for any $t\in \bbc^{*}.$ This induces the action of $\bbc^{*}$ on
$\mathrm{Ext}_A^{1}(X,Y)\setminus\{0\}.$  For any $t\in \bbc^{*},$
we have $t.\varepsilon$ is the following short exact sequence:
$$
\xymatrix{0\ar[r]& Y\ar[rr]^{\left(%
\begin{array}{c}
  1 \\
  0 \\
\end{array}%
\right)}&&L(t.d)\ar[rr]^{\left(%
\begin{array}{cc}
  0 & 1 \\
\end{array}%
\right)}&&X\ar[r]&0}
$$
where $L(t.d)_{\alpha}=\left(%
\begin{array}{cc}
  Y_{\alpha} & td(\alpha) \\
  0 & X_{\alpha} \\
\end{array}%
\right)$  for any $\alpha\in Q_1.$  The orbit space is denoted by
$\mathbb{P}\ext^1_{A}(X,Y)$ and the orbit of $\varepsilon \in
\ext^1_{A}(X,Y)$ is denoted by $\mathbb{P}\varepsilon.$ For a
$G_{\ud}$-invariant constructible subset $\mo$ of $\bbe_{\ud}(A),$
we set $\mathrm{Ext}^{1}_{A}(X,Y)_{\mo}$ to be the subset of
$\mathrm{Ext}^1_{A}(X,Y)$ consisting of the equivalence classes of
extensions with middle terms belonging to $\mo.$

The above $\bbc^{*}$-action on the extensions induces an action on
the middle terms. As a vector space, $L=Y\oplus X.$ So we can
define $t.(y, x)=(ty, x)$ for any $t\in \bbc^{*}$ and $x\in X,
y\in Y$ \cite[Section 5.4]{GLS2006} or \cite[Lemma 1]{Riedtmann}.
For any $L_1\subseteq L,$ this action yields a submodule $t.L_1$
of $L$ isomorphic to $L_1.$ In general, if
$\mathfrak{f}_{L}=(L\supseteq L_1\supseteq L_2\supseteq \cdots
\supseteq L_m=0)$ is a flag of submodules of $L$, then
$t.\mathfrak{f}_{L}=(L\supseteq t.L_1\supseteq
t.L_2\supseteq\cdots\supseteq t.L_m=0)$. Hence, we obtain an
action of $\bbc^{*}$ on the flag of $L.$

\subsection{The product of two evaluation forms}
Let $A_{\ud}:=A_{\ud}(\mathcal{S})$ be the constructible subset of
$\bbe_{\ud}(A)$ as in Section 1.1. For any module $M\in
\bbe_{\ud}(A),$ let $Gr_{\ue}(M)$ be the subvariety of
$Gr_{\ue}(\bbc^{\ud}):=\prod_{i\in Q_0}Gr_{e_i}(\bbc^{d_i})$
consisting of submodules of $M$ with dimension vector
$\ue=(e_i)_{i\in Q_0},$ and let $Gr_{\ue}(\bbe_{\ud}(A))$ be the
constructible subset of $\bbe_{\ud}(A)\times Gr_{\ue}(\bbc^{\ud})$
consisting of pairs $(M,M_1)$ such that $M_1\in Gr_{\ue}(M).$
\begin{Prop}\label{finite}
Let $\ud$ and $\ue$ be two dimension vectors. Then the function
$gr(\ue, \ud): \bbe_{\ud}(A)\rightarrow \bbc$ sending $M$ to
$\chi(Gr_{\ue}(M))$ is a $G_{\ud}$-invariant constructible
function.
\end{Prop}
\begin{proof}Consider the projection: $\phi:Gr_{\ue}(\bbe_{\ud}(A))\rightarrow \bbe_{\ud}(A)$ mapping $(M, M_1)$ to $M.$
It is clear that $\phi$ is a morphism of varieties. By Theorem
\ref{Joyce}, $gr(\ue,\ud)=\phi_{*}(1_{Gr_{\ue}(\bbe_{\ud}(A))})$
is constructible.
\end{proof}
For fixed $\ud,$  we can make finitely many choices of $\ue$ such
that $Gr_{\ue}(\bbe_{\ud}(A))$ is nonempty. This implies the
following corollary.
\begin{Cor}\label{partition}
There is a finite subset $S(\ud)$ of $A_{\ud}$ such that
$A_{\ud}=\bigcup_{i\in S(\ud)}\mo(\ud)_i$ where all $\mo(\ud)_i$
are constructible subsets of $A_{\ud}$ satisfying that for any
$M,M'\in \mo(\ud)_i,$ $\chi(Gr_{\ue}(M))=\chi(Gr_{\ue}(M'))$ for
any $\ue.$
\end{Cor}

 Let $\mathcal{M}(\ud)$ be the vector space over $\bbc$ spanned
by the constructible functions $gr(\ue, \ud)$ for any dimension
vector $\ue.$ For any $M\in A_{\ud},$ we define the evaluation
form $\delta_M: \mathcal{M}(\ud)\rightarrow \bbc$ which maps the
constructible function $gr(\ue)$ to
$\chi(Gr_{\ue}(M))=gr(\ue)(M).$ Using the notations in
\cite{GLS2006}, we set $\str{L}:=\mo(\ud)_i$ for arbitrary $L\in
\mo(\ud)_i.$ Indeed, $\delta_{L}=\delta_{L'}$ for any $L,L'\in
\mo(\ud)_i.$ By abuse of notation, we have $A_{\ud}=\bigcup_{L\in
S(\ud)}\str{L}.$

Let $M, N$ be $A$-modules and $\ue_1, \ue_2$ be dimension vectors.
Fix $M_1\in Gr_{\ue_1}(M), N_1\in Gr_{\ue_2}(N),$ we consider the
natural map:
$$
\beta_{N_1,M_1}: \ext^{1}_{A}(N, M_1)\rightarrow
\ext^{1}_{A}(N,M)\oplus \ext^{1}_{A}(N_1,M_1)
$$
mapping $\varepsilon_{*}\in \ext^{1}_{A}(N, M_1)$ to
$(\varepsilon,\varepsilon')$ such that the following diagram
commutes:
$$
\xymatrix{ \varepsilon': &
0\ar[r]&M_1\ar[r]\ar@{=}[d]&L''\ar[r]\ar[d]&N_1\ar[r]\ar[d]&0\\
\varepsilon_{*}:&0\ar[r]&M_1 \ar[r]\ar[d]&L'\ar[r]\ar[d]&N\ar[r]\ar@{=}[d]&0\\
\varepsilon: & 0\ar[r]&M\ar^-{i}[r]&L\ar[r]^-{\pi}&N\ar[r]&0 }
$$
where $L$ and $L''$ are the pushout and pullback, respectively.
Define
$$
EF^g_{\ue_1,\ue_2}(N,M)=\{(M_1,N_1,\varepsilon, L_1)\mid M_1\in
Gr_{\ue_1}(M), N_1\in Gr_{\ue_2}(N),$$$$\varepsilon\neq 0\in
\ext^1_{A}(N,M)_L\cap Im\beta_{N_1,M_1}, L_1\in
Gr_{\ue_1+\ue_2}(L),L_1\cap i(M)=i(M_1),\pi(L_1)=N_1\}
$$
and
$EF^g_{\ue}(N,M)=\bigcup_{\ue_1+\ue_2=\ue}EF^g_{\ue_1,\ue_2}(N,M).$
By the discussion in Section \ref{C*}, the action of $\bbc^{*}$ on
$\ext^{1}_{A}(N,M)\setminus\{0\}$ naturally induces the action on
$EF^g_{\ue}(N,M)$ by setting
$$
t.(M_1,N_1,\varepsilon, L_1)=(M_1,N_1,t.\varepsilon, t.L_1)
$$
for $(M_1,N_1,\varepsilon, L_1)\in EF^g_{\ue}(N,M)$ and $t\in
\bbc^{*}.$  We denote its orbit space by
$\mathbb{P}EF^g_{\ue}(N,M).$ We also set the evaluation form
$\delta: \mathcal{M}\rightarrow \bbc$ mapping $gr(\ue)$ to
$\chi(\mathbb{P}EF^g_{\ue}(N,M)).$
\begin{theorem}\label{formula1}
Let $M, N\in\mc(\mathcal{S})$. We have
$$
\chi(\mathbb{P}\mathrm{Ext}^{1}_{A}(M,N))\delta_{M\oplus
N}=\sum_{L\in
S(\ud)}\chi(\mathbb{P}\mathrm{Ext}^{1}_{A}(M,N)_{\str{L}})\delta_{L}+\delta.
$$
\end{theorem}
\begin{proof}
Since (for example, see \cite{CC} or
\cite{DXX})$$\chi(Gr_{\ue}(M\oplus
N))=\sum_{\ue_1+\ue_2=\ue}\chi(Gr_{\ue_1}(M))\cdot
\chi(Gr_{\ue_2}(N)),
$$the above formula has the following reformulation
$$
\hspace{-4cm}\chi(\mathbb{P}\ext^{1}_{A}(M,N))\sum_{\ue_1+\ue_2=\ue}\chi(Gr_{\ue_1}(M))\cdot
\chi(Gr_{\ue_2}(N))$$$$\hspace{2.5cm}=\sum_{L\in
S(\ud)}\chi(\mathbb{P}\ext^{1}_{A}(M,N)_{\str{L}})\chi(Gr_{\ue}(L))+\chi(\mathbb{P}EF^g_{\ue}(N,M)).
$$

Now we prove the above reformulation. Define
$$
EF(M,N)=\{(\varepsilon, L_1)\mid \varepsilon\in
\ext^1_{A}(M,N)_{L}\setminus\{0\}, L_1\in Gr_{\ue}(L) \}.
$$
The action of $\bbc^{*}$ on $\ext^{1}_{A}(M,N)$ naturally induces
the action on $EF(M,N)$ \cite[section 5.4 ]{GLS2006}. Under the
action of $\bbc^{*},$ it has the geometric quotient:
$$
\pi:EF(M,N)\rightarrow \mathbb{P}EF(M,N).
$$
We have the natural projection:
$$
p:\mathbb{P}EF(M,N)\rightarrow \mathbb{P}\ext^1_{A}(M,N).
$$
Using Proposition \ref{Euler}, we have $$
\chi(\mathbb{P}EF(M,N))=\sum_{L\in
S(\ud)}\chi(\mathbb{P}\ext^1_{A}(M,N)_{\str{L}})\chi(Gr_{\ue}(L)).
$$
Given $(\varepsilon, L_1)\in EF(M,N),$ let $\varepsilon$ be the
equivalence class of the following short exact sequence:
$$
\xymatrix{\varepsilon:\quad 0\ar[r]& N\ar[rr]^{\left(%
\begin{array}{c}
  1 \\
  0 \\
\end{array}%
\right)}&&L\ar[rr]^{\left(%
\begin{array}{cc}
  0 & 1 \\
\end{array}%
\right)}&&M\ar[r]&0}.
$$
As a vector space, $L=N\oplus M$ and $L_1$ is the subspace of $L.$
We put $M_1=(0, 1)(L_1)$ and $N_1=(1, 0)(L_1)$. It is clear that
$M_1$ and $N_1$ are the submodules of $M$ and $N$, respectively.
Then there is a natural morphism
$$
\phi_0: EF(M,N)\rightarrow
\bigcup_{\ue_1+\ue_2=\ue}Gr_{\ue_1}(M)\times Gr_{\ue_2}(N)
$$
defined by mapping $(\varepsilon, L_1)$ to $(M_1, N_1).$
Furthermore, we have $$\phi_0((\varepsilon,
L_1))=\phi_0(t.(\varepsilon, L_1))$$ for any $(\varepsilon,
L_1)\in EF(M,N)$ and $t\in \bbc^*.$ This induces the morphism
$$
\phi: \mathbb{P}EF(M,N)\rightarrow
\bigcup_{\ue_1+\ue_2=\ue}Gr_{\ue_1}(M)\times Gr_{\ue_2}(N).
$$
 Now we computer the fibre of this morphism
for $M_1\in Gr_{\ue_1}(M)$ and $N_1\in Gr_{\ue_2}(N).$ Consider
the following linear map dual to $\beta_{M_1,N_1}.$
$$
\beta'_{M_1,N_1}: \ext^{1}_{A}(M,N)\oplus
\ext^{1}_{A}(M_1,N_1)\rightarrow \ext^{1}(M_1,N)
$$
mapping $(\varepsilon,\varepsilon')$ to
$\varepsilon_{M_1}-\varepsilon'_{N}$ where $\varepsilon_{M_1}$ and
$\varepsilon'_{N}$ are induced by the inclusions $M_1\subseteq M$
and $N_1\subseteq N,$ respectively as follows:
$$
\xymatrix{
\varepsilon_{M_1}:&0\ar[r]&N\ar[r]\ar@{=}[d]&L_1\ar[r]\ar[d]&M_1\ar[r]\ar[d]&0\\
\varepsilon: & 0\ar[r]&N\ar[r]&L\ar[r]^-{\pi}&M\ar[r]&0 }
$$
where $L_1$ is the pullback, and
$$
\xymatrix{
\varepsilon':&0\ar[r]&N_1\ar[r]\ar[d]&L'\ar[r]\ar[d]&M_1\ar[r]\ar@{=}[d]&0\\
\varepsilon'_N: & 0\ar[r]&N\ar[r]&L'_1\ar[r]&M_1\ar[r]&0 }
$$
where $L'_1$ is the pushout. It is clear that
$\varepsilon,\varepsilon'$ and $M_1,N_1$ induce the inclusions
$L_1\subseteq L$ and $L'\subseteq L'_1.$  And
$$
p_0: \ext^{1}_{A}(M,N)\oplus \ext^{1}_{A}(M_1,N_1)\rightarrow
\ext^{1}_{A}(M,N)
$$
is a projection. By a similar discussion as \cite[Lemma
2.4.2]{GLS2006}, we know
$$p(\phi^{-1}((M_1,N_1)))=\mathbb{P}(p_0(\mathrm{Ker}(\beta'_{M_1,N_1}))).$$
Moreover, by \cite[Lemma 7]{Hubery2005},  for fixed
$\varepsilon\in p_0(\mathrm{Ker}(\beta'_{M_1,N_1})),$ let
$\mathbb{P}\varepsilon$ be its orbit in
$\mathbb{P}(p_0(\mathrm{Ker}(\beta'_{M_1,N_1}))),$ we have
$$
p^{-1}(\mathbb{P}\varepsilon)\cap \phi^{-1}((M_1,N_1))\cong
\hom(M_1,N/N_1).
$$
Using Proposition \ref{Euler}, we get
$$
\chi(\phi^{-1}((M_1,N_1)))=\chi(\mathbb{P}(p_0(\mathrm{Ker}(\beta'_{M_1,N_1})))=\dim_{\bbc}p_0(\mathrm{Ker}(\beta'_{M_1,N_1})).
$$
In the same way, we consider the projection
$$
\varphi:\mathbb{P}EF^g_{\ue}(N,M)\rightarrow
\bigcup_{\ue_1+\ue_2=\ue}Gr_{\ue_1}(M)\times Gr_{\ue_2}(N).
$$
Then
$$
\chi(\varphi^{-1}((M_1,N_1)))=\mathrm{dim}_{\bbc}\ext^1_{A}(N,M)\cap
\mathrm{Im}\beta_{N_1,M_1}.
$$
Now, depending on the fact that $\mc(\mathcal{S})$ is of
Ext-symmetry, we have
$$
\dim_{\bbc}p_0(\mathrm{Ker}(\beta'_{M_1,N_1}))+\mathrm{dim}_{\bbc}\ext^1_{A}(N,M)\cap
\mathrm{Im}\beta_{N_1,M_1}=\mathrm{dim}_{\bbc}\ext^1_{A}(M,N).
$$
Using Proposition \ref{Euler} again, we complete the proof of the
theorem.
\end{proof}
\section{The multiplication formula}
The formula in the last section is not so `symmetric' as the
Geiss-Leclerc-Schr\"oer formula. In order to overcome this
difficulty, we should consider  flags of composition series
instead of Grassmannians of submodules as in \cite{GLS2006}. In
this section, we prove a multiplication formula as an analog of
the Geiss-Leclerc-Schr\"oer formula in \cite{GLS2006}.

Let $A=\bbc Q/\mathcal{I}$ be an algebra associated to a finite
and connected quiver $Q$ and $\mathcal{S}=\{S_1, \cdots, S_n \}$
be a finite set of finite dimensional simple $A$-modules. Let
$\mc(\mathcal{S})$ be a  full subcategory of Ext-symmetry of
$\mathrm{mod}(A)$ associated to $\mathcal{S}.$

Let $A_{\ud}$ be the constructible subset of $\bbe_{\ud}(A)$
consisting of $A$-modules in $\mc(\mathcal{S})$ with dimension
vector $\ud.$ Let $\mathcal{X}$ be the set of pairs $(\textbf{j},
\textbf{c})$ where $\textbf{c}=(c_1,\cdots,c_m)\in \{0,1\}^m$ and
$\textbf{j} = (j_1,\cdots ,j_m)$ is a sequence of integers such
that $S_{j_k}\in \mathcal{S}$ for $1\leq k\leq m.$ Given  $x\in
A_{\ud}$ and $(\textbf{j}, \textbf{c})\in \mathcal{X}$, we define
a $x$-stable flag of type $(\textbf{j},\textbf{c})$ as a
composition series of $x$
$$
\mathfrak{f}_{x} = \left(V=(\bbc^{\ud},x) \supseteq V^1 \supseteq
\cdots \supseteq V^m = 0\right)
$$
of $A$-submodules of $V$ such that $ |V^{k-1}/V^k| =c_kS_{j_k} $
where $S_{j_k}$ is the simple module in $\mathcal{S}.$ Let
$\Phi_{\textbf{j},\textbf{c},x}$ be the variety of $x$-stable
flags of type $(\textbf{j},\textbf{c}).$ We simply write
$\Phi_{\textbf{j}, x}$ when $\textbf{c}=(1,1,\cdots,1).$ Define
$$\Phi_{\textbf{j}}(A_{\ud})=\{(x,\mathfrak{f})\mid x\in
A_{\ud},\mathfrak{f}\in \Phi_{\textbf{j}, x}\}.$$  As in
Proposition \ref{finite}, we consider a projection: $p:
\Phi_{\textbf{j}}(A_{\ud})\rightarrow A_{\ud},$ the function
$p_{*}(1_{\Phi_{\textbf{j}}(A_{\ud})})$ is constructible by
Theorem \ref{Joyce}.
\begin{Prop}\label{const2}
For any type $\textbf{j},$ the function $A_{\ud}\rightarrow \bbc$
mapping $x$ to $\chi(\Phi_{\textbf{j},x})$ is constructible.
\end{Prop}

Let $d_{\textbf{j},\textbf{c}}: \bbe_{\ud}(A)\rightarrow \bbc$ be
the function defined by
 $d_{\textbf{j},\textbf{c}}(x)=\chi(\Phi_{\textbf{j},\textbf{c},x})$ for $x\in
 \bbe_{\ud}(A).$ It is a constructible function as Proposition
 \ref{const2}.
 We simply write
$d_\textbf{j}$ if $\textbf{c}=(1,\cdots,1).$ Define
$\mathcal{M}(\ud)$ to be the vector space spanned by
$d_{\textbf{j}}.$ For fixed $A_{\ud},$ there are finitely many
types $\textbf{j}$ such that $\Phi_{\textbf{j}}(A_{\ud})$ is not
empty. Hence, there exists a finite subset $S(\ud)$ of $A_{\ud}$
such that
$$A_{\ud}=\bigcup_{M\in S(\ud)}\str{M},$$ where $\str{M}=\{M'\in A_{\ud}\mid \chi(\Phi_{\textbf{j},M'})=\chi(\Phi_{\textbf{j},M}) \mbox{ for any type }\textbf{j}\}.$

For any $M\in A_{\ud},$ we define the evaluation form
$\delta_{M}:\mathcal{M}(\ud)\rightarrow \bbc$ mapping a
constructible function $f\in \mathcal{M}(\ud)$ to $f(M).$ We have
$$
\str{M}=\{M'\in A_{\ud}\mid \delta_{M'}=\delta_M\}.
$$
\begin{lemma}\label{directsum}
For $M, N\in \mc(\mathcal{S})$, we have $\delta_{M\oplus
N}=\delta_{M}\cdot \delta_{N}.$
\end{lemma}
The lemma is equivalent to show that
$$\chi(\Phi_{\textbf{j},M\oplus N})=\sum_{\textbf{c}'+\textbf{c}''\sim
1}\chi(\Phi_{\textbf{j},\textbf{c}',M})\cdot\chi(\Phi_{\textbf{j},\textbf{c}'',N}).
$$
Here, $\textbf{c}'+\textbf{c}''\sim 1$ means that $c'_k+c''_k=1$
for $k=1, \cdots, m.$ The proof of the lemma depends on the fact
that under the action of $\bbc^{*}$, $\Phi_{\textbf{j},M\oplus N}$
and its stable subset have the same Euler characteristic. We refer
to \cite{DXX} for details.

The following formula is just the multiplication formula in
\cite[Theorem 1]{GLS2006} when $A$ is a preprojective algebra and
$\mathcal{S}$ is the set of all simple $A$-modules.
\begin{theorem}\label{formula2}
With the above notation, for $M, N\in \mc(\mathcal{S})$,  we have
$$
\chi(\mathbb{P}\mathrm{Ext}^{1}_{A}(M,N))\delta_{M\oplus
N}=\sum_{L\in
S(\ue)}(\chi(\mathbb{P}\mathrm{Ext}^{1}_{A}(M,N)_{\str{L}})+\chi(\mathbb{P}\mathrm{Ext}^{1}_{A}(N,M)_{\str{L}}))\delta_{L},
$$
where $\ue=\underline{\mathrm{dim}}M+\underline{\mathrm{dim}}N.$
\end{theorem}
In the proof of Theorem \ref{formula1}, a key point is to consider
the linear maps $\beta_{M_1,N_1}$ and $\beta'_{M_1,N_1}$ dual to
each other by the property of Ext-symmetry. Now we extend this
idea to the present situation as in \cite{GLS2006}. Let $$
\mathfrak{f}_{M}= (M = M_0 \supseteq M_1 \supseteq \cdots
\supseteq M_m = 0)
$$
be a flag of type $(\textbf{j},\textbf{c}')$ and let
$$\mathfrak{f}_{N}= (N = N_0 \supseteq N_1 \supseteq \cdots
\supseteq N_m = 0)$$ be a flag of type $(\textbf{j},\textbf{c}'')$
such that $c'_k+c''_k=1$ for $k=1,\cdots,m.$ We write
$\textbf{c}'+\textbf{c}''\sim 1.$ For $k=1,\cdots,m,$ let
$\iota_{M,k}$ and $\iota_{N,k}$ be the inclusion maps $M_k \to
M_{k-1}$ and $N_k \to N_{k-1}$, respectively. Define \cite
[Section 2]{GLS2006}
$$
\beta_{\textbf{j},\textbf{c}',\textbf{c}'',\mathfrak{f}_{M},\mathfrak{f}_{N}}:
\bigoplus_{k=0}^{m-2}\ext^1_{A}(N_k,M_{k+1})\rightarrow
\bigoplus_{k=0}^{m-2}\ext^1_{A}(N_k,M_{k})
$$
by the following map
$$
\xymatrix{N_{k}\ar[r]^-{\varepsilon_k}\ar[d]_-{\iota_{N,k}}\ar[dr]&M_{k+1}[1]\ar[d]^-{\iota_{M,k+1}}\\N_{k-1}\ar[r]^-{\varepsilon_{k-1}}&M_{k}[1]}
$$
satisfying
$$
\beta_{\textbf{j},\textbf{c}',\textbf{c}'',\mathfrak{f}_{M},\mathfrak{f}_{N}}(\varepsilon_0,\cdots,\varepsilon_{m-2})=\iota_{M,1}\circ\varepsilon_0+\sum_{k=1}^{m-2}(\iota_{M,k+1}\circ\varepsilon_k-\varepsilon_{k-1}\circ
\iota_{N,k}).
$$
Depending on the fact that $\mc(\mathcal{S})$ is of Ext-symmetry,
we can write down its dual.
$$
\beta'_{\textbf{j},\textbf{c}',\textbf{c}'',\mathfrak{f}_{M},\mathfrak{f}_{N}}:
\bigoplus_{k=0}^{m-2}\ext^1_{A}(M_k,N_{k})\rightarrow\bigoplus_{k=0}^{m-2}\ext^1_{A}(M_{k+1},N_{k})
$$
by the following map
$$
\xymatrix{M_{k+1}\ar[r]^-{\eta_{k+1}}\ar[d]_-{\iota_{M,k+1}}\ar[dr]&N_{k+1}[1]\ar[d]^-{\iota_{N,k+1}}\\M_{k}\ar[r]^-{\eta_k}&N_{k}[1]}
$$
satisfying
$$
\beta'_{\textbf{j},\textbf{c}',\textbf{c}'',\mathfrak{f}_{M},\mathfrak{f}_{N}}(\eta_0,\cdots,\eta_{m-2})=\sum_{k=0}^{m-3}(\eta_k\circ\iota_{M,k+1}-\iota_{N,k+1}\circ\eta_{k+1}
)+\eta_{m-2}\circ \iota_{M,m-1}.
$$
Now, we prove Theorem \ref{formula2}.
\begin{proof}
Define
$$
EF_{\textbf{j}}(M,N)=\{(\varepsilon,\mathfrak{f})\mid
\varepsilon\in \ext^1_{A}(M,N)_{L}, L\in A_{\ue}, \mathfrak{f}\in
\Phi_{\textbf{j},L}\}.
$$
The action of $\bbc^{*}$ on $\mathrm{Ext}^1_A(M,N)$ induces an
action on $EF_{\textbf{j}}(M,N)$. The orbit space under the action
of $\bbc^{*}$ is denoted by $\mathbb{P}EF_{\textbf{j}}(M,N)$ and
the orbit of $(\varepsilon,\mathfrak{f})$ is denoted by
$\mathbb{P}(\varepsilon,\mathfrak{f}).$ We have the natural
projection
$$
p: \mathbb{P}EF_{\textbf{j}}(M,N)\rightarrow
\mathbb{P}\ext^1_{A}(M,N).
$$
The fibre for any $\mathbb{P}\varepsilon\in
\mathbb{P}\ext^1_{A}(M,N)_{L}$ is isomorphic to
$\Phi_{\textbf{j},L}.$ By Theorem \ref{Euler}, we have
$$
\chi(\mathbb{P}EF_{\textbf{j}}(M,N))=\sum_{L\in
S(\ue)}\chi(\mathbb{P}\ext^1_{A}(M,N)_{\str{L}})\chi(\Phi_{\textbf{j},L}).
$$
We also have the natural morphism
$$
\phi:\mathbb{P}EF_{\textbf{j}}(M,N)\rightarrow
\bigcup_{\textbf{c}'+\textbf{c}''\sim
1}\Phi_{\textbf{j},\textbf{c}',M}\times
\Phi_{\textbf{j},\textbf{c}'',N}
$$
mapping $\mathbb{P}(\varepsilon,\mathfrak{f})$ to
$(\mathfrak{f}_M,\mathfrak{f}_N)$ where
$(\mathfrak{f}_M,\mathfrak{f}_N)$ is naturally induced by
$\varepsilon$ and $\mathfrak{f}$ and
$t.(\varepsilon,\mathfrak{f})$ induces the same
$(\mathfrak{f}_M,\mathfrak{f}_N)$ for any $t\in \bbc^{*}.$ By
\cite[Lemma 2.4.2]{GLS2006}, we know
$$p(\phi^{-1}(\mathfrak{f}_M,\mathfrak{f}_N))=\mathbb{P}(p_0(\mathrm{Ker}(\beta'_{\textbf{j},\textbf{c}',\textbf{c}'',\mathfrak{f}_{M},\mathfrak{f}_{N}})),$$
where $p_0: \bigoplus_{k=0}^{m-2}\ext^1_{A}(M_k,N_{k})\rightarrow
\ext^1_{A}(M,N)$ is a projection. On the other hand, by
\cite[Lemma 7]{Hubery2005}, the morphism
$$
p\mid_{\phi^{-1}(\mathfrak{f}_M,\mathfrak{f}_N)}:
\phi^{-1}(\mathfrak{f}_M,\mathfrak{f}_N)\rightarrow
\mathbb{P}(p_0(\mathrm{Ker}(\beta'_{\textbf{j},\textbf{c}',\textbf{c}'',\mathfrak{f}_{M},\mathfrak{f}_{N}}))
$$
has the fibres isomorphic to an affine space. Hence, by Theorem
\ref{Euler},we have
$$
\chi(\phi^{-1}(\mathfrak{f}_M,\mathfrak{f}_N))=\chi(\mathbb{P}(p_0(\mathrm{Ker}(\beta'_{\textbf{j},\textbf{c}',\textbf{c}'',\mathfrak{f}_{M},\mathfrak{f}_{N}})).
$$

Dually, we define
$$
EF_{\textbf{j}}(N,M)=\{(\varepsilon,\mathfrak{f})\mid
\varepsilon\in \ext^1_{A}(N,M)_{L}, L\in A_{\ue}, \mathfrak{f}\in
\Phi_{\textbf{j},L}\}.
$$
The orbit space under $\bbc^{*}$-action is denoted by
$\mathbb{P}EF_{\textbf{j}}(N,M).$ We have the natural projection
$$
q: \mathbb{P}EF_{\textbf{j}}(N,M)\rightarrow
\mathbb{P}\ext^1_{A}(N,M).
$$
The fibre for any $\mathbb{P}\varepsilon\in
\mathbb{P}\ext^1_{A}(N,M)_{L}$ is isomorphic to
$\Phi_{\textbf{j},L}.$ By Theorem \ref{Euler}, we have
$$
\chi(\mathbb{P}EF_{\textbf{j}}(N,M))=\sum_{L\in
S(\ue)}\chi(\mathbb{P}\ext^1_{A}(N,M)_{\str{L}})\chi(\Phi_{\textbf{j},L}).
$$
As in the proof of Theorem \ref{formula1}, there is a natural
morphism
$$
\varphi_0: EF_{\textbf{j}}(N,M)\rightarrow
\bigcup_{\textbf{c}'+\textbf{c}''\sim
1}\Phi_{\textbf{j},\textbf{c}',M}\times
\Phi_{\textbf{j},\textbf{c}'',N}
$$
such that
$$
\varphi_0((\varepsilon,\mathfrak{f}))=\varphi_0(t.(\varepsilon,\mathfrak{f}))
$$
for any $(\varepsilon,\mathfrak{f})\in EF_{\textbf{j}}(N,M)$ and
$t\in \bbc^{*}.$ Hence, we have the morphism
$$
\varphi:\mathbb{P}EF_{\textbf{j}}(N,M)\rightarrow
\bigcup_{\textbf{c}'+\textbf{c}''\sim
1}\Phi_{\textbf{j},\textbf{c}',M}\times
\Phi_{\textbf{j},\textbf{c}'',N}.
$$
 By \cite[Lemma
2.4.3]{GLS2006}, we know
$$
q(\varphi^{-1}(\mathfrak{f}_M,\mathfrak{f}_N))=\mathbb{P}\ext^{1}_{A}(N,M)\cap
Im(\beta_{\textbf{j},\textbf{c}',\textbf{c}'',\mathfrak{f}_{M},\mathfrak{f}_{N}}).
$$
Similar to the above dual situation, by \cite[Lemma
7]{Hubery2005}, the morphism:
$$
q\mid_{\varphi^{-1}(\mathfrak{f}_M,\mathfrak{f}_N)}:
\varphi^{-1}(\mathfrak{f}_M,\mathfrak{f}_N)\rightarrow
\mathbb{P}\ext^{1}_{A}(N,M)\cap
Im(\beta_{\textbf{j},\textbf{c}',\textbf{c}'',\mathfrak{f}_{M},\mathfrak{f}_{N}})
$$
has the fibres isomorphic to an affine space. Hence, by
Proposition \ref{Euler},we have
$$
\chi(\varphi^{-1}(\mathfrak{f}_M,\mathfrak{f}_N))=\chi(\mathbb{P}\ext^{1}_{A}(N,M)\cap
Im(\beta_{\textbf{j},\textbf{c}',\textbf{c}'',\mathfrak{f}_{M},\mathfrak{f}_{N}})).
$$
However, since
$\beta_{\textbf{j},\textbf{c}',\textbf{c}'',\mathfrak{f}_{M},\mathfrak{f}_{N}}$
and
$\beta'_{\textbf{j},\textbf{c}',\textbf{c}'',\mathfrak{f}_{M},\mathfrak{f}_{N}}$
are dual to each other, we have
$$(p_0(\mathrm{Ker}(\beta'_{\textbf{j},\textbf{c}',\textbf{c}'',\mathfrak{f}_{M},\mathfrak{f}_{N}}))^{\bot}=\ext^{1}_{A}(N,M)\cap
Im(\beta_{\textbf{j},\textbf{c}',\textbf{c}'',\mathfrak{f}_{M},\mathfrak{f}_{N}}).$$
Thus we have
\begin{eqnarray}
   && \chi(\mathbb{P}(p_0(\mathrm{Ker}(\beta'_{\textbf{j},\textbf{c}',\textbf{c}'',\mathfrak{f}_{M},\mathfrak{f}_{N}}))+\chi(\mathbb{P}\ext^{1}_{A}(N,M)\cap
Im(\beta_{\textbf{j},\textbf{c}',\textbf{c}'',\mathfrak{f}_{M},\mathfrak{f}_{N}})) \nonumber\\
  &=& \mathrm{dim}_{\bbc}\ext^{1}_{A}(M,N).\nonumber
\end{eqnarray}
Therefore, using Proposition \ref{Euler}, we obtain
$$
\mathbb{P}EF_{\textbf{j}}(M,N)+\mathbb{P}EF_{\textbf{j}}(N,M)=\mathrm{dim}_{\bbc}\ext^{1}(M,N)\cdot\sum_{\textbf{c}'+\textbf{c}''\sim
1}\chi(\Phi_{\textbf{j},\textbf{c}',M})\cdot\chi(\Phi_{\textbf{j},\textbf{c}'',N}).
$$
Now, we have obtained the identity
$$
\mathrm{dim}_{\bbc}\ext^{1}(M,N)\cdot\sum_{\textbf{c}'+\textbf{c}''\sim
1}\chi(\Phi_{\textbf{j},\textbf{c}',M})\cdot\chi(\Phi_{\textbf{j},\textbf{c}'',N})=
$$
$$
\sum_{L\in
S(\ue)}\chi(\mathbb{P}\ext^1_{A}(M,N)_{\str{L}})\chi(\Phi_{\textbf{j},L})+\sum_{L\in
S(\ue)}\chi(\mathbb{P}\ext^1_{A}(N,M)_{\str{L}})\chi(\Phi_{\textbf{j},L})
$$
for any type $\textbf{j}.$ Using Lemma \ref{directsum} and
Proposition \ref{Euler}, we finish the proof of Theorem
\ref{formula2}.
\end{proof}

\section{Examples}
In this section, we give some examples of module subcategories of
Ext-symmetry.

\nd(I) Let $A$ be a preprojective algebra associated to a
connected quiver $Q$ without loops. Let $\mathcal{S}$ be the set
of all simple $A$-modules. Then $\mc(\mathcal{S})$ is of
Ext-symmetry \cite[Theorem 3]{GLS2006}.

\nd (II) Let $A=\bbc Q/\str{\alpha\alpha^{*}-\alpha^{*}\alpha}$ be
an associative algebra associated to the following quiver
$$
Q:= \xymatrix{\bullet\ar@(lu,ld)_{\alpha}\ar@(ru,rd)^{\alpha^*}}
$$
Let $M=(\bbc^m, X_{\alpha}, X_{\alpha^*})$ and $N=(\bbc^n,
Y_{\alpha}, Y_{\alpha^*})$ be two finite dimensional $A$-modules.
Following the characterization of $\mathrm{Ext}^1_{A}(M,N)$ in
Section \ref{C*}, we consider the following isomorphism between
complexes (see \cite[Lemma 1]{CB} or \cite[Section 8.2]{GLS2006}):
$$
\xymatrix@-0.4pc{
\mathrm{Hom}_{\bbc}(M_{\bullet},N_{\bullet})\ar[r]^-{d^0_{M,N}}\ar[d]^{\id}&
\mathrm{Hom}_{\bbc}(M_{\bullet},N_{\bullet})\bigoplus
\mathrm{Hom}_{\bbc}(M_{\bullet},N_{\bullet})
\ar[r]^-{d^1_{M,N}}\ar[d]^{(\id, -\id)}&
\mathrm{Hom}_{\bbc}(M_{\bullet},N_{\bullet})\ar[d]^{-\id}
\\
\mathrm{Hom}_{\bbc}(M_{\bullet},N_{\bullet})\ar[r]^-{d^{1,*}_{N,M}}&
\mathrm{Hom}_{\bbc}(M_{\bullet},N_{\bullet})\bigoplus\mathrm{Hom}_{\bbc}(M_{\bullet},N_{\bullet})\ar[r]^-{d^{0,*}_{N,M}}&
\mathrm{Hom}_{\bbc}(M_{\bullet},N_{\bullet})}
$$
where $M_{\bullet}=\bbc^m, N_{\bullet}=\bbc^n$. Here, we define
$$
d^0_{M,N}(A)=(Y_{\alpha}A-AX_{\alpha},
Y_{\alpha^*}A-AX_{\alpha^*}), d^1_{M,N}(B,
B^*)=Y_{\alpha^{*}}B+B^*X_{\alpha}-Y_{\alpha}B^*-BX_{\alpha^*},
$$
$$
d^{0,*}_{N,M}(B,B^*)=BX_{\alpha^{*}}+B^*X_{\alpha}-Y_{\alpha^*}B-Y_{\alpha}B^*,
d^{1,*}_{N,M}(A)=(Y_{\alpha}A-AX_{\alpha},
-Y_{\alpha^*}A+AX_{\alpha^*})
$$
for any $n\times m$ matrices $A, B$ and $B^*.$ The second complex
is dual to the complex
$$
\xymatrix@-0.4pc{
\mathrm{Hom}_{\bbc}(N_{\bullet},M_{\bullet})\ar[r]^-{d^{0}_{N,M}}&
\mathrm{Hom}_{\bbc}(N_{\bullet},M_{\bullet})\bigoplus\mathrm{Hom}_{\bbc}(N_{\bullet},M_{\bullet})\ar[r]^-{d^{1}_{N,M}}&
\mathrm{Hom}_{\bbc}(N_{\bullet},M_{\bullet})}
$$
with respect to the non-degenerate bilinear form $$\Phi:
\mathrm{Hom}_{\bbc}(N_{\bullet},M_{\bullet})\times
\mathrm{Hom}_{\bbc}(N_{\bullet},M_{\bullet})\rightarrow \bbc$$
mapping $(X, Y)$ to $tr(XY)$.
 As in Section \ref{C*}, we have functorially
$$\mathrm{Ext}_A^1(M,N)=\mathrm{Ker}(d^1_{M,N})/\mathrm{Im}(d^0_{M,N})
\mbox{ and }
\mathrm{DExt}_A^1(N,M)=\mathrm{Ker}(d^{0,*}_{N,M})/\mathrm{Im}(d^{1,*}_{N,M}).$$
Hence, we have a bifunctorial isomorphism:
$$
\mathrm{Ext}_A^1(M,N)\cong \mathrm{DExt}_A^1(N, M).
$$

\nd (III) Deformed preprojective algebras were introduced by
Crawley-Boevey and Holland in \cite{CBH}. Fix
$\lambda=(\lambda_i)_{i\in Q_0}$ where $\lambda_i\in \bbc$. The
deformed preprojective algebra of weight $\lambda$ is an
associative algebra
$$
A(\lambda)=\bbc \overline{Q}/\str{\sum_{\alpha\in
Q_1}\alpha\alpha^{*}-\alpha^{*}\alpha)-\sum_{i\in
Q_0}\lambda_ie_i},
$$
where $\overline{Q}=Q\cup Q^{*}$ is the double of a quiver $Q$
without loops. Let $M, N$ be finite dimensional $A$-modules. As in
Section \ref{C*}, we know $D(M,N)$ is just the kernel of the
following linear map:
$$
\bigoplus_{\alpha\in
\overline{Q}_1}\mathrm{Hom}_{\bbc}(M_{s(\alpha)},
N_{t(\alpha)})\xrightarrow{d^1_{M,N}} \bigoplus_{i\in
Q_0}\mathrm{Hom}_{\bbc}(M_i, N_i),
$$
where $d^1_{M,N}$ maps $(f_{\alpha})_{\alpha\in \overline{Q}_1}$
to $(g_i)_{i\in Q_0}$ such that
$$
g_i=\sum_{\alpha\in Q_1,
s(\alpha)=i}(N_{\alpha^*}f_{\alpha}+f_{\alpha^*}M_{\alpha})-\sum_{\alpha\in
Q_1, t(\alpha)=i}(N_{\alpha}f_{\alpha^*}+f_{\alpha}M_{\alpha^*}).
$$
In the same way as in \cite[Section 8.2]{GLS2006}, we obtain a
bifunctorial isomorphism
$$
\mathrm{Ext}_A^1(M,N)\cong \mathrm{DExt}_A^1(N, M).
$$

\nd (IV) It is easy to construct examples of module subcategories
of Ext-symmetry over an algebra which is not of Ext-symmetry. Let
$A=\bbc Q/\str{\beta\beta^{*}-\beta^{*}\beta}$ be a quotient
algebra associated to the quiver
$$
Q:=
\xymatrix{1\ar[r]^{\alpha}&2\ar@/^/[r]^{\beta}&3\ar@/^/[l]^{\beta^*}}
$$
Let $S_1, S_2$ and $S_3$ be finite dimensional simple $A$-modules
associated to three vertices, respectively. Since
$\mathrm{dim}_{\bbc}\mathrm{Ext}^1(S_1, S_2)=1$ and
$\mathrm{Ext}^1(S_2, S_1)=0$, $A$ is not an algebra of
Ext-symmetry. However,  for $\mathcal{S}=\{S_1, S_3\}$ or $\{S_2,
S_3\},$ $\mc(\mathcal{S})$ is of Ext-symmetry.

 \vspace{0.5cm}
 {\parindent0cm \bf Acknowledgements.}\, We are grateful to the referee
for many helpful comments. In particular, Section 3 is added
following the comments. Furthermore, the second author would like
to thank the Max-Planck Institute for Mathematics in Bonn for a
three-month research stay in 2008.

\bibliographystyle{amsplain}

\end{document}